\pdfoutput=1
\RequirePackage{ifpdf}
\ifpdf %We are running pdfTeX in pdf mode
\documentclass[pdftex]{sigma}
\else
\documentclass{sigma}
\fi

\numberwithin{equation}{section}
\newtheorem{Theorem}{Theorem}[section]
\newtheorem{Corollary}[Theorem]{Corollary}
\newtheorem{Lemma}[Theorem]{Lemma}
\newtheorem{Proposition}[Theorem]{Proposition}

\begin{document}

\newcommand{\arXivNumber}{1410.4125}

\allowdisplaybreaks

\renewcommand{\PaperNumber}{014}

\FirstPageHeading

\ShortArticleName{Generalized Convolution Roots of Positive Def\/inite Kernels on Complex Spheres}

\ArticleName{Generalized Convolution Roots\\
of Positive Def\/inite Kernels on Complex Spheres}

\Author{Victor S.~BARBOSA and Valdir A.~MENEGATTO}

\AuthorNameForHeading{V.S.~Barbosa and V.A.~Menegatto}

\Address{Instituto de Ci\^encias Matem\'aticas e de Computa\c c\~ao, Universidade de S\~ao Paulo,\\
Caixa Postal 668, 13560-970, S\~ao Carlos - SP, Brasil}
\Email{\href{mailto:victorrsb@gmail.com}{victorrsb@gmail.com}, \href{menegatt@icmc.usp.br}{menegatt@icmc.usp.br}}

\ArticleDates{Received October 16, 2014, in f\/inal form February 10, 2015; Published online February 13, 2015}

\Abstract{Convolution is an important tool in the construction of positive def\/inite kernels on a~manifold.
This contribution provides conditions on an $L^2$-positive def\/inite and zonal kernel on the unit sphere of
$\mathbb{C}^q$ in order that the kernel can be recovered as a~generalized convolution root of an equally positive
def\/inite and zonal kernel.}


\Classification{33C55; 41A35; 41A63; 42A82; 42A85}

\section{Introduction}

In this paper, we deal with a~specif\/ic problem involving the so-called positive def\/inite kernels.
These kernels have importance in the resolution of issues pertaining to many areas of mathematics, including
approximation theory, functional analysis, statistics, etc, presenting themselves in both theory and applications.
The concept of positive def\/initeness usually appears in two formats as we now introduce.

Let~$X$ be a~nonempty set and~$K$ a~kernel on~$X$, that is, a~complex function~$K$ with domain $X \times X$.
The kernel~$K$ is {\em positive definite} if
\begin{gather*}
\sum\limits_{m,n=1}^k c_m\overline{c_n}K(x_m,x_n) \geq 0,
\end{gather*}
for all positive integer~$k$, all subsets $\{x_1, x_2, \ldots, x_k\}$ of~$X$ and any complex numbers $c_1, c_2, \ldots,
c_k$.
As for the second concept of positive def\/initeness, we need to assume two facts:~$X$ is endowed with a~measure~$\mu$ so
that $(X,\mu)$ is a~$\sigma$-f\/inite measure space and the kernel~$K$ needs to be an element of $L^2(X\times X, \mu\times
\mu)$.
As usual, we will speak of the elements in this space as if they were functions, with equality interpreted as
equality~$\mu$-a.e..
A~kernel~$K$ in $L^2(X\times X, \mu\times \mu)$ is {\em $L^2$-positive definite} if
\begin{gather}
\label{inner}
\int_X \left(\int_X K(x,y) f(y)d\mu(y)\right)\overline{f(x)} d\mu(x) \geq 0,
\qquad
f \in L^2(X,\mu).
\end{gather}
If $\mathcal{K}$ is the integral operator on $L^2(X,\mu)$ generated by~$K$, the previous inequality corresponds to
$\langle \mathcal{K}f,f\rangle_2 \geq 0$, $ f \in L^2(X,\mu)$, in which $\langle \cdot, \cdot \rangle_2$ is the usual
inner product of $L^2(X,\mu)$.

In the cases in which the setting allows a~comparison between the two concepts, it is not hard to see that they do not
coincide.
If one can speak of continuity in~$X$, then a~chance of coincidence increases considerably (see~\cite{ferreira} for
additional information on a~possible equivalence).
In this paper, we prefer to deal almost exclusively with $L^2$-positive def\/initeness since the source of the main
problem to be discussed comes from functional analysis.

If an $L^2$-positive def\/inite kernel~$K$ from $L^2(X\times X, \mu\times \mu)$ is available, a~quite general problem is
the determination of a~kernel~$S$ in $L^2(X\times X, \mu\times \mu)$ satisfying the integral relation
\begin{gather}
\label{recovery}
\int_{X} S(x, \xi)S(\xi,y) d\mu(\xi)=K(x,y),
\qquad
x,y \in X.
\end{gather}
If existence is guaranteed, the uniqueness of~$S$ cannot be usually reached unless additional requirements are added
(indeed, if~$S$ is a~solution, then $-S$ is another one).

In~\cite{ferreira}, the existence of a~positive def\/inite solution for this general problem was solved in the case in
which~$X$ is a~metric space, the measure~$\mu$ is strictly positive and~$K$ is continuous.
The solution~$S$ itself was constructed via the square root $\mathcal{K}^{1/2}$ of the integral operator~$\mathcal{K}$
acting on~$L^2(X,\mu)$ and generated by~$K$.
As a~matter of fact, the operator~$\mathcal{K}^{1/2}$ was an integral operator on~$L^2(X,\mu)$ itself and the generating
kernel of $\mathcal{K}^{1/2}$ was a~solution~$S$, equally continuous and positive def\/inite.
The results proved in~\cite{ferreira} agree with those previously surveyed in~\cite{schaback}, the main dif\/ference
between them being the setting considered in each case.

If~$X$ has an enhanced structure, the original kernel~$K$ may have additional and desirable features.
Hence, the problem specializes to the following version: is there a~solution kernel~$S$ matching all the relevant
additional properties the original kernel~$K$ may have?
In some cases, the left-hand side of~\eqref{recovery} generalizes the notion of convolution.
For instance, if~$X$ is a~compact two-point homogeneous manifold and~$\mu$ is the usual ``volume'' measure on~$X$, then
the integral in~\eqref{recovery} is an extension of the concept of spherical convolution described in~\cite{plat,plat1}
when the kernel~$S$ is zonal (isotropic).
In particular, the solution of the proposed problem boils down to the f\/inding of a~convolution root~$S$ for~$K$ having
the same features as~$K$.
This is relevant on itself because (spherical) convolution is known to be an ef\/f\/icient tool in the construction of
positive def\/inite kernels.
We also observe that the case in which~$X$ is a~sphere, is of particular relevancy in statistics and approximation
theory (see~\cite{estrade,schreiner,ziegel}).

Among the many research problems collected in~\cite{gneiting} the readers can f\/ind open questions aligned with the
material covered in the references mentioned in the previous paragraph and with the results to be developed in Sections
2 and 3 ahead.
In those sections, we will chase a~solution for~\eqref{recovery} in the case in which~$X$ is the unit sphere
$\Omega_{2q}$ in $\mathbb{C}^q$,~$\mu$ is the usual Lebesgue measure in $\Omega_{2q}$, and~$K$ is zonal, following the
direct procedure adopted in~\cite{ziegel}, a~method that contemplates a~Fourier analysis perspective.
The procedure itself encompasses these steps: gene\-ra\-lized convolutions preserve zonality; zonal kernels on $\Omega_{2q}$
have specif\/ic expansions that are necessary and suf\/f\/icient for zonality; generalized convolution acts on these
expansions by squaring coef\/f\/icients; the existence of zonal generalized convolution square roots is a~question of
convergence of expansions; starting from general zonal kernels, one step of generalized convolution generates
a~continuous zonal kernel.
In particular, it should become clear to the reader that the paper deals with a~standard argument that features an
operation (generalized convolution), a~symmetry property (zonality), and a~representation (series expansion) that
maintains the symmetry and is invariant under the operation.

We believe the approach developed here can be adapted to many other situations, a~typical example being some compact
two-point homogeneous manifolds.
The setting and the notion of generalized spherical convolution used in the paper, can be found in~\cite{leves}, where
a~multiplier version of the Bernstein inequality on $\Omega_{2q}$ was obtained.
In order to direct the reader for potential applications involving the analysis on $\Omega_{2q}$ and perhaps, to the
specif\/ic material covered here, we mention~\cite{thiru,torre}.

The spherical approach taken in the paper is probably less general than that adopted in the Hilbert--Schmidt theory (that
was the option in~\cite{ferreira,schaback}).
For the convenience of the reader and also for possible comparisons, that more general approach will be sketched in
Section~\ref{Section4}.
However, it is valuable to mention that the technique adopted here is more explicit and allows the study of regularity
properties of the convolution roots and of the kernel.

\section{The spherical case plus zonality}\label{s1}

Let $\sigma_q$ denote the surface measure on the unit sphere $\Omega_{2q}$ of $\mathbb{C}^q$.
Assuming~$K$ is a~zonal and $L^2$-positive def\/inite kernel on $\Omega_{2q}$, we will seek for a~kernel~$S$ in
$L^2(\Omega_{2q}\times \Omega_{2q}, \sigma_q\times \sigma_q)$, of the same type as~$K$, so that
\begin{gather*}
\int_{\Omega_{2q}} S(x, \xi)S(\xi,y) d\sigma_q(\xi)=K(x,y),
\qquad
x,y \in \Omega_{2q}.
\end{gather*}

Classical results from the analysis on the sphere $\Omega_{2q}$ and also basic Fourier analysis on $\Omega_{2q}$ will be
needed along the way.
In more than half of this section, we will detail some of them and will prove some others not available in the
literature.
Recent references related to the material to be described in this section are~\cite{bezu,wunsche} while classical ones
are~\cite{koo1,rudin}.

The measure $\sigma_q$ is invariant with respect to the group $\mathcal{O}(2q)$ of all unitary linear operators on
$\mathbb{C}^q$ in the following sense: if~$E$ is a~$\sigma_q$-measurable subset of $\Omega_{2q}$ and~$\rho$ is an
element of $\mathcal{O}(2q)$ then $\sigma_q(\rho(E))=\sigma_q(E)$.
Zonality of a~kernel~$K$ refers to the invariance property
\begin{gather*}
K(\rho(x),\rho(y))=K(x,y),
\qquad
x,y \in \Omega_{2q},
\qquad
\rho \in \mathcal{O}(2q).
\end{gather*}

We simplify notation by writing $L^2(\Omega_{2q}):=L^2(\Omega_{2q}, \sigma_q)$, $\Omega_{2q}^2:=\Omega_{2q}\times
\Omega_{2q}$ and $L^2(\Omega_{2q}^2):=L^2(\Omega_{2q}^2,\sigma_q \times \sigma_q)$.
The pertinent inner products in these Hilbert spaces will appear in normalized form, that is,
\begin{gather*}
\langle f,g \rangle_2:=\frac{1}{\omega_q}\int_{\Omega_{2q}}f(x)\overline{g(x)} d\sigma_q(x),
\qquad
f,g \in L^2(\Omega_{2q}),
\end{gather*}
and
\begin{gather*}
\langle K_1, K_2 \rangle_2:=\frac{1}{\omega_q^2}\int_{\Omega_{2q}^2}K_1(x,y)\overline{K_2(x,y)} d(\sigma_q\times
\sigma_q)(x,y),
\qquad
K_1,K_2 \in L^2\big(\Omega_{2q}^2\big),
\end{gather*}
respectively, in which $\omega_q:=2\pi^q/{(q-1)!}$ is the surface area of $\Omega_{2q}$.

The generalized convolution of two kernels $K_1$ and $K_2$ from $L^2(\Omega_{2q}^2)$ is the kernel $K_1 * K_2$ given~by
the formula
\begin{gather}
\label{convol}
(K_1 * K_2)(x,y)=\frac{1}{\omega_q}\int_{\Omega_{2q}}K_1(x, \xi)K_2(\xi,y) d\sigma_{q}(\xi),
\qquad
x,y \in \Omega_{2q}.
\end{gather}
Clearly, the def\/inition makes sense as long as the integral is well-def\/ined.

From now on, in order to make the presentation clearer, we need to consider two separate cases.
For $q\geq 2$, the zonality of~$K$ corresponds to the existence of a~function $K':B[0,1] \to \mathbb{C}$ so that
\begin{gather}
\label{zonal}
K(x,y)=K'(x \cdot y),
\qquad
x,y \in \Omega_{2q},
\end{gather}
in which $\cdot$ denotes the usual inner product in $\mathbb{C}^q$ and $B[0,1]:=\{z \in \mathbb{C}: z\overline{z}\leq
1\}$.
Thus, if $K_2$ is zonal and we replace $K_1$ with a~function~$f$ that depends upon $\xi$ only, then the convolution $K_1
* K_2$ becomes the spherical convolution of $K_2'$ and $\overline{f}$ as def\/ined in~\cite{oliveira,quinto}.

Adapting arguments found in~\cite{menegatto}, one can see that a~kernel~$K$ on $\Omega_{2q}$ is $L^2$-positive def\/inite
and zonal if and only if the generating function $K'$ appearing in~\eqref{zonal} have a~double series representation of
the form
\begin{gather}
\label{coeff}
K'(z)=\sum\limits_{m,n=0}^\infty a_{m,n}^{q-2}(K') R_{m,n}^{q-2}(z),
\qquad
z \in B[0,1],
\end{gather}
in which $a_{m,n}^{q-2}(K') \geq 0$, $m,n \in \mathbb{Z}_+$, with convergence in $L^2(B[0,1],\nu_{q-2})$, where
\begin{gather*}
d\nu_{q-2}(z)=\frac{q-1}{\pi}\big(1-|z|^2\big)^{q-2} dxdy,
\qquad
z \in B[0,1].
\end{gather*}
The symbol $R_{m,n}^{q-2}$ stands for the disk or generalized Zernike polynomial of bi-degree $(m,n)$ associated with
the dimension~$q$, that is,
\begin{gather*}
R_{m,n}^{q-2}(z):=r^{|m-n|}e^{i(m-n)\theta}R_{m\wedge n}^{(q-2,|m-n|)}\big(2r^2 -1\big),
\qquad
z=re^{i\theta}\in B[0,1],
\end{gather*}
where $R_{m\wedge n}^{(q-2,|m-n|)}$ is the Jacobi polynomial of degree $m\wedge n:=\min\{m,n\}$ associated with the pair
of numbers $q-2$ and $|m-n|$, normalized so that $R_{m\wedge n}^{(q-2,|m-n|)}(1)=1$.
Since
\begin{gather*}
R_{k}^{(q-2,|m-n|)}=\frac{P_{k}^{(q-2,|m-n|)}}{P_{k}^{(q-2,|m-n|)}(1)},
\end{gather*}
in which $P_k^{(q-2,|m-n|)}$ is the classical Jacobi polynomial of degree~$k$ associated with $q-2$ and $|m-n|$, as
def\/ined in~\cite{szego}, then $|R_{m,n}^{q-2}(z)|\leq 1$, $z \in B[0,1]$ and $R_{m,n}^{q-2}(1)=1$ for all~$m$ and~$n$.
An alternative formula for $R_{m,n}^{q-2}$ via standard monomials is
\begin{gather*}
R_{m,n}^{q-2}(z)=\frac{m!n!(q-2)!}{(m+q-2)!(n+q-2)!}\sum\limits_{j=0}^{m\wedge n}\frac{(-1)^j
(m+n+q-2-j)!}{j!(m-j)!(n-j)!}z^{m-j}\overline{z}^{n-j}.
\end{gather*}
The set $\{R_{m,n}^{q-2}:m,n \in \mathbb{Z}_+\}$ is an orthogonal system in $L^2(B[0,1],\nu_{q-2})$, that is,
\begin{gather}
\label{orthog}
\int_{B[0,1]}R_{m,n}^{q-2}(z)\overline{R_{k,l}^{q-2}(z)} d\nu_{q-2}(z)=\frac{\delta_{mk}\delta_{nl}}{h_{m,n}^{q-2}},
\end{gather}
where
\begin{gather*}
h_{m,n}^{q-2}=\frac{m+n+q-1}{q-1}\binom{m+q-2}{q-2}\binom{n + q-2}{q-2}.
\end{gather*}
If~$K$ is a~zonal kernel on $\Omega_{2q}$, we will write $a_{m,n}^{q-2}(K')$ to denote the Fourier coef\/f\/icient of the
gene\-ra\-ting function $K'$, with respect to the orthogonal basis $\{R_{m,n}^{q-2}: m,n \in \mathbb{Z}_+\}$ of
$L^2(B[0,1],\nu_{q-2})$:
\begin{gather*}
%\label{foueq}
a_{m,n}^{q-2}(K')=h_{m,n}^{q-2}\int_{B[0,1]}K'(z)\overline{R_{m,n}^{q-2}(z)} d\nu_{q-2}(z),
\qquad
m,n \in \mathbb{Z}_+.
\end{gather*}
The complex Funk--Hecke formula in $\Omega_{2q}$ establishes an intimate connection among everything we have mentioned so
far.
If~$Y$ is a~spherical harmonic of bi-degree $(m,n)$ in~$q$ dimensions, that is,~$Y$ is the restriction to $\Omega_{2q}$
of a~polynomial of bidegree $(m,n)$ in $\mathbb{C}^q$ that belongs to the kernel of the complex Laplacian in
$\mathbb{C}^q$, and $K'$ is a~function in $L^2(B[0,1],\nu_{q-2})$, it states that \cite{oliveira,quinto}
\begin{gather*}
\frac{1}{\omega_q}\int_{\Omega_{2q}}K'(x \cdot y) \overline{Y(x)} d\sigma_q(x)
=\bigg[\int_{B[0,1]}K'(z)\overline{R_{m,n}^{q-2}(z)} d\nu_{q-2}(z)\bigg] \overline{Y(y)},\ y \in \Omega_{2q},
\end{gather*}
as long as the integral on the left-hand side exists.
In particular, the Funk--Hecke formula simplif\/ies to
\begin{gather*}
\frac{1}{\omega_q}\int_{\Omega_{2q}}K'(x \cdot y) \overline{Y(x)} d\sigma_q(x)=\frac{a_{m,n}^{q-2}(K')}{h_{m,n}^{q-2}}
\overline{Y(y)},
\qquad
y \in \Omega_{2q}.
\end{gather*}

We now move to the case $q=1$.
The zonality of a~kernel~$K$ on $\Omega_2$ corresponds to the existence of a~function $K': \Omega_2 \to \mathbb{C}$ so that
\begin{gather}
\label{zonal2}
K(x,y)=K'(x \cdot y),
\qquad
x,y \in \Omega_{2},
\end{gather}
The characterization of the $L^2$-positive def\/inite and zonal kernels on $\Omega_2$ skips a~little bit from the previous
representation in higher dimensions.
In this particular case, the function $K'$ appearing in~\eqref{zonal2} have a~series representation of the form
\begin{gather*}
%\label{coeffi}
K'(z)=\sum\limits_{k=-\infty}^\infty a_{k}(K') z^{k},
\qquad
z \in \Omega_2,
\end{gather*}
in which $a_{k}(K') \geq 0$, $k\in \mathbb{Z}$, with convergence of the series in $L^2(\Omega_2)$.
If one wants to speak of orthogonality, integration needs to be done in accordance with the standard Fourier theory,
interpreting $\Omega_2$ as the quotient space $\mathbb{R}/2\pi\mathbb{Z}$, under the equivalence relation $x \sim y$ if
and only if $x-y \in 2\pi\mathbb{Z}$.
If $q:\mathbb{R} \to \Omega_2$ is the quotient map, then (the invariant) integration on $\Omega_{2}$ reads like
\begin{gather*}
\int_{\Omega_{2}}f(x)d\sigma_2(x):=\int_{0}^{2\pi} f(q(\theta))d\theta,
\end{gather*}
in which $d\theta$ is the Lebesgue measure element on the interval $[0,2\pi)$.
In this sense, if
\begin{gather*}
R_m(z):=z^m,
\qquad
z \in \Omega_2,
\qquad
m\in \mathbb{Z},
\end{gather*}
then $\{R_m: m \in \mathbb{Z}\}$ is an orthonormal basis of $L^2(\Omega_{2q},\sigma_2)$.
Indeed, f\/irst observe that
\begin{gather*}
\frac{1}{\omega_2}\int_{\Omega_2} R_m(z)\overline{R_n(z)} d\sigma_2(z)=\frac{1}{\omega_2}\int_0^{2\pi}
e^{i(m-n)\theta}d\theta.
\end{gather*}
As so, if $m=n$ then the integral on the right-hand side is just 1.
Otherwise, we can pick $\theta_0 \in [0,2\pi]$ so that $e^{i(m-n)\theta_0} \neq 1$ and use a~change of variables to
justify the equalities
\begin{gather*}
\frac{1}{\omega_2}\int_0^{2\pi} e^{i(m-n)\theta}d\theta=\frac{1}{\omega_2}\int_0^{2\pi}
e^{i(m-n)(\theta+\theta_0)}d\theta=\frac{1}{\omega_2}e^{i(m-n)\theta_0}\int_0^{2\pi} e^{i(m-n)\theta}d\theta,
\end{gather*}
and arrive at
\begin{gather*}
\frac{1}{\omega_2}\int_0^{2\pi} e^{i(m-n)\theta}d\theta=0.
\end{gather*}

To close the section, we will list some specif\/ic properties of the space
\begin{gather*}
\mathcal{B}^2\big(\Omega_{2q}^2\big):=\big\{K \in L^2\big(\Omega_{2q}^2\big): K \; \text{is zonal}\big\},
\end{gather*}
some of them to be used ahead.
The norm in the spaces $L^2(\Omega_{2q}^2)$ will be written as $\|\cdot\|_2$.
The f\/irst one is of general interest and is included for completeness only.
If $\{K_n\}$ is a~sequence in $\mathcal{B}^2(\Omega_{2q}^2)$ converging to $K \in L^2(\Omega_{2q}^2)$, $\rho \in
\mathcal{O}(2q)$ and we write $\rho K$ to denote the kernel $(x,y) \in \Omega_{2q} \times \Omega_{2q} \to K(\rho x, \rho
y)$, the zonality of each $K_n$ plus the invariance of $\sigma_{q}$ with respect to~$\rho$ justify the inequality
\begin{gather*}
0\leq \|K-\rho K\|_2   \leq   \|K-K_n\|_2 +\|K_n - \rho K\|_2
\\
\phantom{0}
=\|K-K_n\|_2 +\|\rho K_n - \rho K\|_2=2\|K-K_n\|_2.
\end{gather*}
Letting $n \to \infty$ leads to $K=\rho K$ in $L^2(\Omega_{2q}^2)$ and, consequently, we can conclude that
$\mathcal{B}^2(\Omega_{2q}^2)$ is a~closed subspace of $L^2(\Omega_{2q}^2)$.

To proceed, we will need two additional notations:
\begin{gather*}
Z_{m,n}(x,y):=R_{m,n}^{q-2}(x\cdot y),
\qquad
m,n \in \mathbb{Z}_+,
\qquad
x,y \in \Omega_{2q},
\qquad
q\geq 2,
\end{gather*}
and
\begin{gather*}
Z_m(x,y):=R_m(x\cdot y)=x^my^{-m},
\qquad
m\in \mathbb{Z},
\qquad
x,y \in \Omega_{2}.
\end{gather*}
Since~$q$ will remain f\/ixed, the omission of the dimension~$q$ in both notations introduced above should cause no
confusion.

In the proposition below, we establish orthogonality properties of the sets $\{Z_{m,n}: m,n \in \mathbb{Z}_+\}$ and
$\{Z_m: m\in \mathbb{Z}\}$ in the Hilbert space $(\mathcal{B}^2(\Omega_{2q}^2), \langle \cdot, \cdot \rangle_2)$.

\begin{Proposition}
\label{orthogo}
The following assertions hold:
\begin{enumerate}\itemsep=0pt
\item[$(i)$] $(q\geq 2)$ the Hilbert space $\mathcal{B}^2(\Omega_{2q}^2)$ is isometrically isomorphic to
$L^2(B[0,1],\nu_{q-2})$;

\item[$(ii)$] $(q\geq 2)$ the set $\{(h_{m,n}^{q-2})^{1/2}Z_{m,n}: m,n \in \mathbb{Z}_+\}$ is an orthonormal basis of
$\mathcal{B}^2(\Omega_{2q}^2)$;

\item[$(iii)$] the set $\{Z_m: m\in \mathbb{Z}\}$ is an orthonormal basis of $\mathcal{B}^2(\Omega_{2}^2)$.
\end{enumerate}
\end{Proposition}

\begin{proof}
Let us consider the case $q\geq 2$ f\/irst.
If $K \in \mathcal{B}^2(\Omega_{2q}^2)$, the Funk--Hecke formula implies that
\begin{gather*}
\|K\|_2^2=\frac{1}{\omega_q^2}\int_{\Omega_{2q}^2}K'(x\cdot y)\overline{K'(x\cdot y)}d(\sigma_q \times
\sigma_q)(x,y)
\\
\phantom{\|K\|_2^2}
=\int_{B[0,1]}K'(z)\overline{K'(z)}d\nu_{q-2}(z)=\|K'\|_{2,q}^2,
\end{gather*}
in which $\|\cdot \|_{2,q}$ stands for the usual norm in the space $L^2(B[0,1],\nu_{q-2})$.
In particular, $K\in \mathcal{B}^2(\Omega_{2q}^2) \to K' \in L^2(B[0,1],\nu_{q-2})$ is an isometry.
This takes care of (i).
Since $\{R_{m,n}^{q-2}: m,n\in \mathbb{Z}_+\}$ is an orthogonal basis in $L^2(B[0,1],\nu_{q-2})$, the polarization
identity reveals that its inverse image by the isometry in $(i)$ is an orthogonal basis in
$\mathcal{B}^2(\Omega_{2q}^2)$.
Recalling~\eqref{orthog}, the orthonormality of the set in (ii) follows.
In order to see it is complete in $\mathcal{B}^2(\Omega_{2q}^2)$, let $K \in \mathcal{B}^2(\Omega_{2q}^2)$ and assume
that
\begin{gather*}
\frac{1}{\omega_q^2}\int_{\Omega_{2q}^2}K(x,y)\overline{R_{m,n}(x,y)} d(\sigma_q \times \sigma_q)(x,y)=0,
\qquad
m,n \in \mathbb{Z}_+.
\end{gather*}
Since that corresponds to
\begin{gather*}
\frac{1}{\omega_q}\int_{\Omega_{2q}}\bigg(\frac{1}{\omega_q}\int_{\Omega_{2q}}K'(x \cdot y)R_{n,m}^{q-2}(x\cdot
y)d\sigma_q(x)\bigg)d\sigma_q(y)=0,
\qquad
m,n \in \mathbb{Z}_+,
\end{gather*}
an application of the Funk--Hecke formula leads to the equality
\begin{gather*}
\big(h_{m,n}^{q-2}\big)^{-1}a_{m,n}^{q-2}(K')=0,
\qquad
m,n \in \mathbb{Z}_+,
\end{gather*}
that is,
\begin{gather*}
\int_{B[0,1]}K'(z)\overline{R_{m,n}^{q-2}(z)} d\nu_{q-2}(z)=0,
\qquad
m,n \in \mathbb{Z}_+.
\end{gather*}
Since $\{R_{m,n}^{q-2}:m,n \in \mathbb{Z}_+\}$ is an orthogonal basis of $L^2(B[0,1],\nu_{q-2})$ and this space is
complete, then $K'=0$ in $L^2(B[0,1],\nu_{q-2})$.
In particular, $K=0$ in $\mathcal{B}^2(\Omega_{2q}^2)$ and (ii) is proved.
The orthonormality assertion in (iii) is clear.
As for completeness, let $K \in \mathcal{B}^2(\Omega_{2}^2)$ and assume that
\begin{gather*}
\frac{1}{\omega_2^2}\int_{\Omega_{2}^2}K(x,y)\overline{Z_{m}(x,y)} d(\sigma_2 \times \sigma_2)(x,y)=0,
\qquad
m \in \mathbb{Z}.
\end{gather*}
Since this equality is precisely
\begin{gather*}
\frac{1}{4\pi^2}\int_{\Omega_{2}}\left(\int_{\Omega_{2}}K'(x \cdot y)x^{-m}d\sigma_2(x)\right)y^md\sigma_2(y)=0,
\qquad
m,n \in \mathbb{Z},
\end{gather*}
it follows that
\begin{gather*}
a_m(K')\int_{\Omega_2}y^{-m} y^m d\sigma_2(y)=0,
\qquad
m \in \mathbb{Z},
\end{gather*}
in which $a_m(K')$ is the~$m$-th usual Fourier coef\/f\/icient of $K'$.
In other words, $a_m(K')=0$, $m \in \mathbb{Z}$, that is, $K'=0$.
\end{proof}

We close the section detaching an obvious consequence of Proposition~\ref{orthogo}(ii).

\begin{Lemma}[$q\geq 2$]\label{end3}
If $K'$ belongs to $L^2(B[0,1],\nu_{q-2})$, then $\sum\limits_{m,n=0}^{\infty}
\big(h_{m,n}^{q-2}\big)^{-1}a_{m,n}^{q-2}(K')^2$ is convergent.
\end{Lemma}

\section[Convolution roots in $\mathcal{B}^2(\Omega_{2q}^2)$]{Convolution roots in $\boldsymbol{\mathcal{B}^2(\Omega_{2q}^2)}$}\label{Section3}

In this section, we f\/inally analyze the problem described in the f\/irst paragraph of the previous section.
The f\/irst steps provide technical results on the convolution operation $*$ def\/ined in~\eqref{convol}, in the case when
one or both kernels are the basic elements of Proposition~\ref{orthogo}.
The main results of the paper come right after that.

\begin{Lemma}
\label{selfa}
The following properties hold:
\begin{enumerate}\itemsep=0pt
\item[$(i)$] the convolution of a~hermitian kernel from $L^2(\Omega_{2q}^2)$ with itself is $L^2$-positive definite on~$\Omega_{2q}$;
\item[$(ii)$] if $K_1$ and $K_2$ belong to $\mathcal{B}^2(\Omega_{2q}^2)$, then $K_1 * K_2$ does so.
\end{enumerate}
\end{Lemma}
\begin{proof}
Write $\mathcal{K} * \mathcal{K}$ to denote the integral operator on $L^2(\Omega_{2q})$ generated by $K * K$.
If $f \in L^2(\Omega_{2q})$ and~$K$ is a~kernel on $\Omega_{2q}$, then a~double application of Fubini's theorem implies
that
\begin{gather*}
\langle \mathcal{K}*\mathcal{K}(f),f\rangle_2=\frac{1}{\omega_q^3}\int_{\Omega_{2q}}\int_{\Omega_{2q}}K(x,
\xi)\overline{f(x)} d\sigma_q(x) \int_{\Omega_{2q}} K(\xi,y) f(y)d\sigma_q(y) d\sigma_{q}(\xi).
\end{gather*}
If~$K$ is hermitian, the previous equality reduces itself to
\begin{gather*}
\langle
\mathcal{K}*\mathcal{K}(f),f\rangle_2=\frac{1}{\omega_q}\int_{\Omega_{2q}}\bigg|\frac{1}{\omega_q}\int_{\Omega_{2q}}K(x,
\xi)\overline{f(x)} d\sigma_q(x)\bigg|^2 d\sigma_{q}(\xi).
\end{gather*}
This takes care of~(i).
As for (ii), if $\rho \in \mathcal{O}_{2q}$, then
\begin{gather*}
(K_1 * K_2)(\rho (x),\rho (y))=\frac{1}{\omega_q}\int_{\Omega_{2q}}K_1(x, \rho^* \xi)K_2(\rho^* \xi,y) d\sigma_{q}(\xi),
\qquad
x,y \in \Omega_{2q},
\end{gather*}
in which $\rho^*$ is the adjoint of~$\rho$.
Since the measure $\sigma_q$ is invariant with respect to elements of~$\mathcal{O}_{2q}$, it follows that
\begin{gather*}
(K_1 * K_2)(\rho x,\rho y)=\frac{1}{\omega_q}\!\int_{\Omega_{2q}}\!\!\! K_1(x, \rho^* \xi)K_2(\rho^* \xi,y)d\sigma_{q}(\rho^* \xi)
=(K_1*K_2)(x,y),
\qquad\!\!\!
x,y \in \Omega_{2q},
\end{gather*}
and~$K$ is zonal.
\end{proof}

The lemma below intends to produce formulas for convolutions involving the basis elements of $\mathcal{B}^2(\Omega_{2q}^2)$.

\begin{Lemma}
\label{lema*}
The following formulas hold:
\begin{gather*}
h_{k,l}^{q-2} (Z_{m,n} * Z_{k,l})=\delta_{km}\delta_{ln}Z_{k,l},
\qquad
m,n,k,l \in \mathbb{Z}_+,
\\
Z_m * Z_n=\delta_{m,n} Z_m,
\qquad
m,n \in \mathbb{Z}.
\end{gather*}
\end{Lemma}
\begin{proof}
In the case $q\geq 2$, an application of the Funk--Hecke formula leads to
\begin{gather*}
(Z_{m,n} * Z_{k,l})(x,y)=\frac{1}{\omega_q}\int_{\Omega_{2q}}R_{m,n}^{q-2}(x \cdot \xi)R_{k,l}^{q-2}(\xi \cdot y)d\sigma_{q}(\xi)
\\
\phantom{(Z_{m,n} * Z_{k,l})(x,y)}
=\frac{1}{\omega_q} \int_{\Omega_{2q}}R_{k,l}^{q-2}(\xi \cdot y) \overline{R_{m,n}^{q-2}(\xi \cdot x)}d\sigma_{q}(\xi)
\\
\phantom{(Z_{m,n} * Z_{k,l})(x,y)}
=\frac{a_{m,n}^{q-2}\big(R_{k,l}^{q-2}\big)}{h_{m,n}^{q-2}}\overline{R_{m,n}^{q-2}(y \cdot x)},
\qquad
x,y \in \Omega_{2q}.
\end{gather*}
It is now clear that
\begin{gather*}
(Z_{m,n} * Z_{k,l})(x,y)=\frac{\delta_{mk}\delta_{nl}}{h_{m,n}^{q-2}}R_{m,n}^{q-2}(x\cdot y),
\qquad
x,y \in \Omega_{2q}.
\end{gather*}
Moving to the case $q=1$, direct computation shows that
\begin{gather*}
(Z_m * Z_n)(x,y)=\frac{1}{\omega_2} \int_{\Omega_2} Z_m(x,\xi)Z_n(\xi,y)d\sigma_2(\xi)
\\
\phantom{(Z_m * Z_n)(x,y)}
=x^my^{-n}\frac{1}{\omega_2}\int_{\Omega_2}\overline{R_m(\xi)}R_n(\xi)\sigma_2(\xi),
\qquad
m,n \in \mathbb{Z}_+,
\qquad
x,y\in \Omega_2.
\end{gather*}
The use of the orthonormality of $\{R_m:m\in \mathbb{Z}\}$ in $L^2(\Omega_2)$ in the previous equation implies the
second formula in the statement of the lemma.
\end{proof}

\begin{Theorem}
\label{k*smn}
If~$K$ belongs to $\mathcal{B}^2(\Omega_{2q}^2)$ then
\begin{gather*}
K*Z_{m,n}=\frac{a_{m,n}^{q-2}(K')}{h_{m,n}^{q-2}}Z_{m,n}=Z_{m,n}*K,
\qquad
m,n \in \mathbb{Z}_+,
\qquad
q\geq 2,
\end{gather*}
and
\begin{gather*}
K*Z_{m}=a_m(K')Z_{m}=Z_{m}*K,
\qquad
m,n \in \mathbb{Z}_+,
\qquad
q=1.
\end{gather*}
\end{Theorem}
\begin{proof}
In the case $q\geq 2$, both equalities follow from calculations via the Funk--Hecke formula.
As for the case $q=1$, it can be done by direct calculations and a~change of variables.
\end{proof}

The formulae in the theorem below are the key step towards the main results of this section.

\begin{Theorem}
\label{k*k}
If~$K$ belongs to $\mathcal{B}^2(\Omega_{2q}^2)$ then
\begin{gather*}
\langle K*K,Z_{m,n}\rangle_2=\bigg[\frac{a_{m,n}^{q-2}(K')}{h_{m,n}^{q-2}}\bigg]^2,
\qquad
m,n \in \mathbb{Z}_+,
\qquad
q\geq 2,
\end{gather*}
and
\begin{gather*}
\langle K*K,Z_{m}\rangle_2=a_m(K')^2,
\qquad
m \in \mathbb{Z},
\qquad
q=1.
\end{gather*}
\end{Theorem}
\begin{proof}
Fix $m,n \in \mathbb{Z}_+$.
Introducing the equality from Lemma~\ref{lema*}, we see that
\begin{gather*}
\langle K*K, Z_{m,n}\rangle_2 =\frac{1}{\omega_q^2}\int_{\Omega_{2q}}\int_{\Omega_{2q}}(K*K)(\xi,\eta)Z_{m,n}(\eta,\xi)d\sigma_q(\xi) d\sigma_q(\eta)
\\
\phantom{\langle K*K, Z_{m,n}\rangle_2}
=\frac{h_{m,n}^{q-2}}{\omega_q^2}\int_{\Omega_{2q}}\int_{\Omega_{2q}}(K*K)(\xi,\eta)(Z_{m,n}
*Z_{m,n})(\eta,\xi)d\sigma_q(\xi) d\sigma_q(\eta).
\end{gather*}
Appealing to the def\/inition of convolution and using Fubini's theorem twice to interchange integration orders, we deduce
that
\begin{gather*}
\langle K*K,
Z_{m,n}\rangle_2=\frac{h_{m,n}^{q-2}}{\omega_q^2}\int_{\Omega_{2q}}\int_{\Omega_{2q}}(K*Z_{m,n})(x,y)(Z_{m,n}*K)(y,x)d\sigma_q(x)\sigma_q(y).
\end{gather*}
Applying Theorem~\ref{k*smn}, we obtain
\begin{gather*}
\langle K*K, Z_{m,n}\rangle_2=\big(h_{m,n}^{q-2}\big)^{-1} a_{m,n}^{q-2}(K')^2 \|Z_{m,n}\|_2^2=\big(h_{m,n}^{q-2}\big)^{-1}
a_{m,n}^{q-2}(K')^2 \big\|R_{m,n}^{q-2}\big\|_{2,q}^2,
\end{gather*}
while~\eqref{orthog} yields the f\/irst equality in the statement of the theorem.
The other equality is proved in a~similar manner.
\end{proof}

\begin{Theorem}
\label{convergen}
Let~$K$ belong to $L^2(\Omega_{2q}^2)$.
If $K=J*J$ for some~$J$ in $\mathcal{B}^2(\Omega_{2q}^2)$, then

\begin{gather*}
\sum\limits_{m,n=0}^{\infty}h_{m,n}^{q-2}|\langle K,Z_{m,n}\rangle_2|<\infty,
\qquad
q\geq 2,
\end{gather*}
and
\begin{gather*}
\sum\limits_{m=-\infty}^{\infty}|\langle K, Z_m \rangle_2|<\infty,
\qquad
q=1.
\end{gather*}
\end{Theorem}
\begin{proof}
If $K=J*J$ with $J \in \mathcal{B}^2(\Omega_{2q}^2)$, then Lemma~\ref{selfa} implies that $K \in
\mathcal{B}^2(\Omega_{2q}^2)$ as well.
In the case $q\geq 2$, that allows an application of Theorem~\ref{k*k} to deduce the equality
\begin{gather*}
\langle K,Z_{m,n}\rangle_2=\bigg[\frac{a_{m,n}^{q-2}(J')}{h_{m,n}^{q-2}}\bigg]^2,
\qquad
m,n \in \mathbb{Z}_+,
\end{gather*}
that is,
\begin{gather*}
h_{m,n}^{q-2}\langle K,Z_{m,n}\rangle_2=\bigg[\frac{a_{m,n}^{q-2}(J')}{(h_{m,n}^{q-2})^{1/2}}\bigg]^2,
\qquad
m,n \in \mathbb{Z}_+.
\end{gather*}
Since $\{(h_{m,n}^{q-2})^{1/2}R_{m,n}^{q-2}: m,n \in \mathbb{Z}^+\}$ is an orthonormal basis of $L^2(B[0,1],\nu_{q-2})$
and $J' \in L^2(B[0,1],\nu_{q-2})$, we invoke Lemma~\ref{end3} to conclude that
\begin{gather*}
\sum\limits_{m,n=0}^{\infty}h_{m,n}^{q-2}|\langle K,Z_{m,n}\rangle_2| <\infty.
\end{gather*}
The case $q=1$ can be handled in a~similar manner.
\end{proof}

Needless to say that, due to Lemma~\ref{selfa}, the modulus sign in the convergence outcome of the previous theorem can
be removed, as long as the kernel~$J$ is hermitian.

Theorem~\ref{convergen} suggests a~criterion for the existence of convolution roots, therefore, a~solution to the
problem described at the beginning of Section~\ref{s1}.

\begin{Theorem}
\label{raizsuf}
Let~$K$ be a~kernel in $L^2(\Omega_{2q}^2)$.
If $q\geq 2$, assume that all the Fourier coefficients $\langle K,Z_{m,n}\rangle_2$ are nonnegative and that
\begin{gather*}
\sum\limits_{m,n=0}^{\infty}h_{m,n}^{q-2}\langle K,Z_{m,n}\rangle_2 <\infty.
\end{gather*}
Otherwise, assume that all Fourier coefficients $\langle K, Z_m \rangle_2$ are nonnegative and that
\begin{gather*}
\sum\limits_{m=-\infty}^{\infty}\langle K, Z_m \rangle_2 <\infty.
\end{gather*}
Then, there exists an $L^2$-positive definite kernel~$P$ in $\mathcal{B}^2(\Omega_{2q}^2)$ such that $K=P * P$.
In particular,~$K$ is an $L^2$-positive definite element of $\mathcal{B}^2(\Omega_{2q}^2)$.
\end{Theorem}
\begin{proof}
We prove the theorem in the case $q\geq 2$ only.
Let us consider the zonal kernel~$P$ for which
\begin{gather*}
P'=\sum\limits_{m,n=0}^{\infty}h_{m,n}^{q-2} \langle K,Z_{m,n}\rangle_2^{1/2}Z_{m,n}.
\end{gather*}
The expansion of $P'$ with respect to the orthonormal basis $\{(h_{m,n}^{q-2})^{1/2} Z_{m,n}:m,n\in \mathbb{Z}_+\}$ of~$\mathcal{B}^2(\Omega_{2q}^2)$ is
\begin{gather*}
P'=\sum\limits_{m,n=0}^{\infty}\big(h_{m,n}^{q-2}\big)^{1/4}\langle
K,\big(h_{m,n}^{q-2}\big)^{1/2}Z_{m,n}\rangle_2^{1/2}\big(h_{m,n}^{q-2}\big)^{1/2}Z_{m,n},
\qquad
m,n \in \mathbb{Z}_+.
\end{gather*}
In particular, our convergence assumption implies that $P'$ belongs to $\mathcal{B}^2(\Omega_{2q}^2)$.
Our additional assumptions and the characterization provided in~\eqref{coeff} implies that~$P$ is $L^2$-positive
def\/inite.
At last, a~help of Theorem~\ref{k*k} leads to
\begin{gather*}
\langle P * P, Z_{m,n}\rangle_2=\bigg[\frac{a_{m,n}^{q-2}(P')}{h_{m,n}^{q-2}}\bigg]^2=\langle K,Z_{m,n}\rangle_2,
\qquad
m,n \in \mathbb{Z}_+.
\end{gather*}
This suf\/f\/ices to conclude that $P*P=K$ in $L^2(\Omega_{2q}^2)$.
\end{proof}

The next result reveals that, in the setting we have adopted, the existence of a~zonal convolution root of the kernel
implies the existence of a~continuous one.

\begin{Proposition}
Let~$K$ belong to $L^2(\Omega_{2q}^2)$.
If $K=J*J$ for some~$J$ in $\mathcal{B}^2(\Omega_{2q}^2)$, then~$K$ is in fact a~continuous kernel.
\end{Proposition}
\begin{proof}
We know already that $K \in \mathcal{B}^2(\Omega_{2q}^2)$.
Hence, it can be written as
\begin{gather}
\label{kserie}
K(x,y)=\sum\limits_{m,n=0}^{\infty} a_{m,n}^{q-2}(K') R_{m,n}^{q-2} (x\cdot y).
\end{gather}
Since
$
a_{m,n}^{q-2}=h_{m,n}^{q-2}\langle K,Z_{m,n} \rangle_2$,
$m,n \in \mathbb{Z}_+$,
it follows that
\begin{gather*}
\sum\limits_{m,n=0}^{\infty}\big|a_{m,n}^{q-2}(K') R_{m,n}^{q-2}(x\cdot y)\big|
\leq \sum\limits_{m,n=0}^{\infty}h_{m,n}^{q-2}|\langle K,Z_{m,n} \rangle_2|,
\qquad
x,y \in \Omega_{2q}.
\end{gather*}
Recalling Theorem~\ref{convergen}, we can apply the Weierstrass M-test to conclude that the series in~\eqref{kserie}
converges uniformly in $\Omega_{2q}^2$.
Since each $R_{m,n}$ is continuous, it follows that the same series def\/ines a~continuous kernel.
\end{proof}

The next result is a~consequence of the previous proposition and Theorem 2.3 in~\cite{ferreira}.

\begin{Corollary}
Let~$K$ belong to $L^2(\Omega_{2q}^2)$.
If $K=J*J$ for some~$J$ in $\mathcal{B}^2(\Omega_{2q}^2)$, then~$K$ is positive definite in the usual sense.
\end{Corollary}

This a~simplif\/ied version of Theorem~\ref{raizsuf}.

\begin{Theorem}
If~$K$ is a~continuous, zonal and $L^2$-positive definite kernel on $\Omega_{2q}$, then there exists an $L^2$-positive
definite kernel~$P$ in $\mathcal{B}^2(\Omega_{2q}^2)$ such that $K=P*P$.
\end{Theorem}
\begin{proof}
We only prove the assertion in the case $q\geq 2$.
Write $K(x,y)=K'(x \cdot y)$, $x,y \in \Omega_{2q}$, with $K'$ as described in~\eqref{coeff}.
Applying the Funk--Hecke formula, we see that
\begin{gather*}
\langle K,Z_{m,n}\rangle_2=\frac{1}{\omega_q^2}\int_{\Omega_{2q}^2}K'(x\cdot y)\overline{R_{m,n}^{q-2}(x\cdot y)}
d(\sigma_q\times \sigma_q)(x,y)
\\
\phantom{\langle K,Z_{m,n}\rangle_2}
=\frac{1}{\omega_q}\frac{a_{m,n}^{q-2}(K')}{h_{m,n}^{q-2}}\int_{\Omega_{2q}}\overline{R_{m,n}^{q-2}(y \cdot y)}
d\sigma_q(y),
\qquad
m,n \in \mathbb{Z}_+,
\end{gather*}
that is,
\begin{gather*}
\langle K,Z_{m,n}\rangle_2=\frac{a_{m,n}^{q-2}(K')}{h_{m,n}^{q-2}}\geq 0,
\qquad
m,n \in \mathbb{Z}_+.
\end{gather*}
Furthermore,
\begin{gather*}
\sum\limits_{m,n=0}^{\infty}h_{m,n}^{q-2} \langle K,Z_{m,n}\rangle_2=\sum\limits_{m,n=0}^{\infty}
a_{m,n}^{q-2}(K')<\infty.
\end{gather*}
An application of Theorem~\ref{raizsuf} justif\/ies the assertion of the theorem.
\end{proof}

We advise the reader that the results obtained in this section can be reproduced in some other important settings in
both, real and complex versions.
As a~matter of fact, the approach can be developed as long as an $L^2$ structure similar to the one used here is
available in the setting to be considered.

\section{Final remarks: the solution via integral operators}\label{Section4}

Here, for the sake of completeness, we recall the functional analysis approach via Hilbert--Schmidt operator usually used
to solve a~version of the general recovery problem~\eqref{recovery}.
One needs the separability of $L^2(X,\mu)$, the continuity of~$K$ and, in addition, both concepts of positive
def\/initeness need to coincide.
In this case, $\eqref{inner}$ holds and the integral operator $\mathcal{K}$ generated by~$K$ is a~Hilbert--Schmidt
operator.

The classical spectral and Mercer's theories guarantee the features below.
\begin{enumerate}\itemsep=0pt
\item[(i)] The operator $\mathcal{K}$ is self-adjoint and Hilbert--Schmidt.

\item[(ii)] There exists a~f\/inite or countably inf\/inite set of positive real numbers $\lambda_{1}(\mathcal{K}) \geq
\lambda_{2}(\mathcal{K}) \geq \dots \geq 0$ and a~corresponding orthonormal system $\{\phi_{n}\}$ of $L^{2}(X,\mu)$ so
that
\begin{gather*}
\mathcal{K}(\phi_n)=\lambda_n(\mathcal{K}) \phi_n,
\qquad
n=1,2,\ldots.
\end{gather*}
\item[(iii)] The sequence $\{\lambda_n(\mathcal{K})\}$ is square-summable and converges to 0.

\item[(iv)] The integral operator $\mathcal{K}$ has an $L^2$ representation in the form
\begin{gather*}
\mathcal{K}(f)=\sum\limits_{n=1}^{\infty}\lambda_n(\mathcal{K})\langle f,\phi_n \rangle_2 \phi_n,
\qquad
f \in L^2(X).
\end{gather*}
\item[(v)] The kernel~$K$ has an expansion in the form
\begin{gather*}
K(x,y)=\sum\limits_{n=1}^{\infty}\lambda_{n}(\mathcal{K})\phi_{n}(x)\overline{\phi_{n}(y)},
\end{gather*}
with convergence in $L^2(X\times X, \mu \times \mu)$.
\end{enumerate}

All this information being available, the following classical result holds.

\begin{Theorem}
\label{main1}
If~$K$ is a~continuous $L^2$-positive definite kernel on~$X$, then there exists an $L^2$-positive definite kernel~$S$ in
$L^2(X\times X, \mu\times \mu)$, which can be taken continuous, so that
\begin{gather*}
\int_{X} S(x, \xi)S(\xi,y) d\mu(\xi)=K(x,y),
\qquad
x,y \in X.
\end{gather*}
\end{Theorem}

The following result is also pertinent.

\begin{Corollary}
Under the setting and the conditions stated in Theorem{\rm ~\ref{main1}}, the square root~$\mathcal{K}^{1/2}$ of~$\mathcal{K}$
is the integral operator on~$L^2(X,\mu)$ generated by~$S$.
\end{Corollary}

In this general context, the kernels~$K$ and~$S$ end up having the same $L^2$ expansion structure.
Hence, it is expectable that~$S$ will have the very same properties~$K$ has, as long as the properties are attached to
the Hilbert space structure of $L^2(X,\mu)$.
In a~certain sense, the results in Section~\ref{Section3} ratify this in the spherical setting, at least when the property to be preserved is zonality.
We are unaware of specif\/ic papers addressing the analysis of a~recovery property similar to the one considered here on~$\Omega_{2q}$ via this Functional Analysis approach.

\subsection*{Acknowledgements}
The f\/irst author was partially supported by CAPES.
The second one by FAPESP, under grant 2014/00277-5.
Special thanks goes to the anonymous referees for the careful reading of the paper and for pointing corrections and
suggestions that led to this f\/inal form of the paper.

\pdfbookmark[1]{References}{ref}
\LastPageEnding

\end{document}